\documentclass[11pt]{article}
\usepackage[cp1251]{inputenc} 
\usepackage[russian]{babel}  
\usepackage[dvips]{graphicx}
\usepackage[dvips]{epsfig}
\usepackage{amssymb}
 \makeindex

\begin{document}

\begin{center}
А.Р. МИРОТИН, М.А. РОМАНОВА
\end{center}

\begin{center}
{\bf  ИНТЕРПОЛЯЦИОННЫЕ МНОЖЕСТВА АЛГЕБРЫ ОБОБЩЕННЫХ АНАЛИТИЧЕСКИХ ФУНКЦИЙ}
\end{center}

\begin{center}
 E-mail: amirotin@yandex.ru
\end{center}

 Исследуются свойства компактных интерполяционных множеств алгебр обобщенных аналитических функций и даются
удобные достаточные условия интерполяционности.

Abstract. The properties of the compactness of interpolation sets of algebras of generalized analytic functions are investigated and convenient sufficient conditions for interpolation are given.

\vspace{5mm}

Известная теорема Рудина-Карлесона (см., напр., \cite{Rud74}, с. 125,
6.3.1) характеризует компактные подмножества окружности, интерполяционные
для диск-алгебры, как множества  нулевой меры Лебега. Однако уже в случае поликруга ситуация
более сложна; для компакта на $n$-мерном торе ($n>1$) наличие нулевой
лебеговой меры (меры Хаара), продолжая быть необходимым, перестает быть достаточным
условием интерполяционности (см. там же,теорема
6.3.4).  Ряд достаточных условий,   смысл которых в малости
рассматриваемого компакта в том или ином смысле, для поликруга изложен в (\cite{Rud74},
параграфы 6.2, 6.3;
\cite{Val}).  Случай шара рассмотрен в \cite{Rud84}, где имеются
дальнейшие библиографические указания.

Интересное обобщение понятия аналитичности, тесно связанное с абстрактным
гармоническим анализом, было предложено Р.~Аренсом и И.~М.~Зингером в работе \cite{AS},
которая породила обширную литературу. При этом основное внимание уделялось
случаю архимедовой упорядоченности (см.\cite{Gam}, глава VII; \cite{Rud62}, глава VIII).
С современным состоянием  этой
области исследований можно познакомиться по монографии \cite{Ton} и
обзору \cite{GT1} (см. также \cite{Gr}).

Целью настоящей работы является исследование  подмножеств компактных
абелевых групп, интерполяционных для алгебры обобщенных аналитических  по
Аренсу-Зингеру функций. Строго говоря, понятие обобщенной
аналитичности, принятое в данной работе,
является несколько менее ограничительным, чем в \cite{AS}, и следует
подходу, предложенному в \cite{Mir}.  Попутно мы устанавливаем некоторые
свойства обобщенных аналитических функций, имеющие, возможно, самостоятельный интерес.
Архимедова упорядоченность, как правило, не предполагается.

Структура статьи следующая. В первом параграфе, который носит
вспомогательный характер, собраны необходимые сведения о полухарактерах
полугрупп. Во втором параграфе доказана теорема о среднем для обобщенных
аналитических функций, и рассматривается свойство полиномиальной аппроксимации
алгебры таких функций.  Основные результаты   содержатся в третьем
параграфе, где исследуются свойства компактных интерполяционных множеств и даются
удобные достаточные условия интерполяционности.
   Результаты статьи были частично анонсированы в  \cite{MR}. Статья продолжает цикл работ автора, посвященных гармоническому анализу на полугруппах,
см. \cite{Mir2}, \cite{Hilb}, \cite{SbMath}.

\vspace{10mm}

\centerline{\bf 1 Полухарактеры}
\vspace{5mm}

В этом параграфе  приведены  сведения о полухарактерах абелевых полугрупп,
используемые в дальнейшем.

Всюду ниже $S$ --- дискретная абелева полугруппа с сокращениями и единицей $e$,
записываемая мультипликативно, $G=S^{-1}S$ --- группа частных для $S$.

    {\it Полухарактером} полугруппы $S$ называется
гомоморфизм $\psi$ полугруппы $S$ в мультипликативную полугруппу
$\overline {\Bbb{D}}=\{ z\in~\Bbb{C}: |z|\leq 1\}$, не являющийся тождественным
нулем.  {\it Характерами}  называются полухарактеры, равные по модулю единице.

    Множество всех полухарактеров далее обозначается $\widehat {S}$, а его подмножество,
состоящее из неотрицательных полухарактеров, --- $\widehat {S}_+$.
 Наделенные топологией поточечной сходимости это компактные топологические
полугруппы по умножению с единицей $1$ ($\widehat {S}$ компактно, например,  как
  пространство максимальных идеалов алгебры $l_1(S)$ \cite{AS}).

Компактную группу всех характеров полугруппы $S$ будем обозначать $X$. Каждый
характер $\chi\in X$ единственным образом продолжается до характера группы
$G$ по формуле $\chi(a^{-1}b):=\overline{\chi(a)}\chi(b)\quad(a,b\in S$),
что позволяет отождествлять $X$ с группой характеров группы $G$ (это
специальный случай результата, установленного в \cite{Mir2}).
Аналогично каждый полухарактер $\rho\in \widehat {S}_+$   единственным образом
продолжается с  $S(\rho)$ до  гомоморфизма $\tilde{\rho}$
группы частных  $G(\rho)$ полугруппы $S(\rho)$    в мультипликативную
полугруппу положительных действительных чисел.

Согласно теореме 3.1 из \cite{AS} каждый полухарактер $\psi\in\widehat {S}$
может быть представлен в виде $\psi=\rho\chi$, где $\chi\in X$, а полухарактер
$\rho\in \widehat {S}_+$ определяется по $\psi$  однозначно.

С  полухарактером $\rho\in \widehat {S}_+$ связаны подполугруппы
$S(\rho):=\{s\in S: \rho(s)> 0\}$ и $S^{\rho}:=\{s\in S: \rho(s)=1\}$
полугруппы $S$, дополнение которых, если оно не пусто, является идеалом
$S$. Такие идеалы  называются {\it простыми};  это в точности те идеалы,
дополнения которых являются подполугруппами полугруппы $S$, поскольку  индикаторы
(характеристические функции) дополнений простых идеалов принадлежат
$\widehat {S}_+$. Простые идеалы,
отличные от $S\setminus\{e\}$, будем называть {\it нетривиальными}.
Через $\omega$ обозначим индикатор множества $\{e\}$.
Наконец отметим,что степень  $\rho^0$ по определению есть
индикатор $S(\rho)$, и что  $\rho^z\in \widehat {S}\setminus X$ при
$\rho\in \widehat {S}_+, \quad \rho\ne 1,\quad z\in \Pi$, где
$\Pi:=\{ {\rm Re} z>0\}$ (\cite{AS},\S 7).

Конец доказательства  будет  обозначаться знаком $\Box$.

\vspace{15mm}
\centerline{\bf 2 Алгебра обобщенных аналитических функций}
\vspace{5mm}

{\bf Определение 1.} Комплекснозначную функцию $F$ на $\widehat {S}
\setminus X$ будем называть {\it обобщенной аналитической},
если для любых полухарактеров $\rho, \psi$ из $\widehat {S}\setminus X$,
$\rho \geq 0$ отображение
$z \mapsto F(\rho^{z} \psi)$
  аналитично в открытой правой полуплоскости $\Pi$ и непрерывно в +0.

Равномерную алгебру функций, непрерывных на $\widehat {S}$ и аналитических на
$\widehat {S}\setminus X$, будем обозначать $A(\widehat {S})$. Для
аддитивной полугруппы $\Bbb{Z}_+^n$ получаем поликруговую алгебру, ---
равномерную алгебру функций,
аналитических в открытом поликруге и непрерывных в его замыкании (в этом
случае определение 1 гарантирует аналитичность по каждой переменной в отдельности).

Каждая аналитическая по Аренсу-Зингеру функция принадлежит $A(\widehat
{S})$ (это можно доказать как в \cite{AS}, теорема 7.4). Следующий пример показывает,
что обратное  утверждение, вообще говоря, неверно.

{\bf Пример 1}. Если $S$ есть аддитивная полугруппа $\{0, 2, 3, \ldots
\}$, то функция $F(z)=z$ принадлежит $A(\widehat {S})=A(\overline
{\Bbb{D}})$, но по теореме 2.6 из \cite{AS} не аналитична по Аренсу-Зингеру,
так как преобразование
Фурье ее сужения на единичную окружность не соcредоточено на $S$ .

Примерами функций из $A(\widehat {S})$ являются {\it аналитические полиномы}, т. е.
 функции на $\widehat {S}$ вида
$$
p(\psi):= \sum_{i=1}^{n}c_{i} \widehat {a_{i}}(\psi),
$$
где $c_{i}\in \Bbb{C}, \psi \in \widehat {S}$,
$\widehat {a_{i}}(\psi):=\psi(a_{i}), a_{i}\in S$.

Более общим образом, алгебре  $A(\widehat {S})$ принадлежит {\it
преобразование Лапласа} ($=$ преобразование Гельфанда \cite{AS}) любой функции $f$ из
$l_1(S)$, т. е. $A(\widehat {S})$ содержит все функции  вида
$$
\widehat {f}(\psi)=\sum\limits_{s\in
S}f(s)\psi(s)\quad(\psi\in\widehat{S}).
$$

В частности, зафиксируем весовую функцию $w:S\to \Bbb{R}_+$,
удовлетворяющую условию $\sum_{s\in S}w(s)=1$. Тогда для любой
ограниченной функции $b$ на $S$ функция $\psi\mapsto\langle\psi,b\rangle\quad$, где
$\psi\in \widehat {S}$,  принадлежит $A(\widehat {S})$
(угловые скобки обозначают скалярное произведение в $l_2(S,w)$; мы неоднократно будем
использовать эти обозначения в дальнейшем).

{\bf Теорема 1}. {\it Пусть $F\in A(\widehat {S})$. Если полухарактеры
$\rho_1,\rho_2\in \widehat {S}_+$ удовлетворяют либо условию $S(\rho_1)=
S(\rho_2)$, либо условию  $S^{\rho_1}=S^{\rho_2}$, то справедливо равенство
$$
\int\limits_X F(\rho_1\chi)d\chi=\int\limits_X F(\rho_2\chi)d\chi,
$$
\noindent
где $d\chi$--- нормированная мера Хаара группы} $X$.

{\bf Доказательство}. Предположим сначала, что $S(\rho_1)= S(\rho_2)$. Тогда
функция
$$
f(z)= \int\limits_X F({\rho_1}^z\rho_2\chi)d\chi
$$
\noindent
будет, очевидно, непрерывна  в  полуплоскости  $\Pi=\{{\rm Re} z>0\}$,
а также в $+0$. Кроме
того, она аналитична в $\Pi$, например, по теореме
Мореры. В силу условий Коши-Римана имеем для любого $\theta>0\quad (z=x+iy)$
$$
f^{'}(\theta)=\frac{\partial f(\theta)}{\partial
x}=\frac{1}{i}\frac{\partial f(\theta)}{\partial y}=
$$
$$
\frac{1}{i}\lim\limits_{h\to 0}\frac{1}{h} (\int\limits_X
F({\rho_1}^{\theta+ih}\rho_2\chi)d\chi -\int\limits_X
F({\rho_1}^\theta\rho_2\chi)d\chi ).
$$
\noindent
Заметим теперь, что характер $\rho_1^{ih}$  может быть
продолжен с полугруппы  $S(\rho_1)$ до некоторого характера $\chi_h$ полугруппы $S$
(достаточно
применить теорему Понтрягина к продолжению характера  $\rho_1^{ih}$   на
группу частных полугруппы $S(\rho_1)$, а затем сузить полученный характер группы $G$ на
$S$ ).   Поскольку
$\rho_1^{\theta+ih}(s)=\rho_1^{\theta}(s)\chi_h(s)$ при всех $s\in S$, то
$f^{'}(\theta)=0$ для любого $\theta>0$ в силу инвариантности меры Хаара.
Так как $\rho_1^{0}\rho_2=\rho_2$, то  с учетом непрерывности $f$ в $+0$
 имеем
$$
0=f(1)-f(0)=\int\limits_X F(\rho_1\rho_2\chi)d\chi -\int\limits_X
F(\rho_2\chi)d\chi.
$$
\noindent
Симмерия между $\rho_1$ и $\rho_2$  завершает первую часть доказательства.

Предположим, наконец, что  $S^{\rho_1}=S^{\rho_2}$. Обозначая через
$\omega_1$ индикатор этой подполугруппы  и применяя результат первой
части к полухарактерам  $\rho_i$ и $\rho_i^n$, имеем в силу теоремы
Лебега о мажорированной сходимости $(i=1,2)$
$$
\int\limits_X F(\rho_i\chi)d\chi=\int\limits_X F({\rho_i}^n\chi)d\chi  \to
\int\limits_X F(\omega_1\chi)d\chi\quad(n\to\infty).   \Box
$$

Из предыдущей теоремы легко вытекает следующий вариант теоремы о среднем для
функций из  $A(\widehat {S})$. Напомним, что через $\omega$ мы
обозначаем индикатор одноточечного множества $\{e\}$,  который принадлежит
$\widehat {S}_+$, если (и только если) $S^{-1}\cap S=\{e\}$.

{\bf Теорема 2.} {\it Пусть  $S^{-1}\cap S=\{e\}$. Если существует такой
полухарактер $\rho_1\in\widehat {S}_+$, что $0<\rho_1(s)<1$ при всех $s\in S, s\ne
e$, то при всех}  $F\in A(\widehat {S})$
$$
\int\limits_X F(\chi)d\chi=F(\omega).
$$

{\bf Доказательство}. Заметим, что $S(1)=S(\rho_1), S^{\rho_1}=S^\omega.$ Следовательно,

$$
\int\limits_X F(1\cdot\chi)d\chi=\int\limits_X F(\rho_1\chi)d\chi=
\int\limits_X F(\omega\chi)d\chi =F(\omega).  \Box
$$

{\bf Определение 2.} Будем говорить, что алгебра $A(\widehat {S})$ {\it обладает свойством
полиномиальной аппроксимации}, если произвольная функция $F$ из $A(\widehat {S})$
может быть равномерно приближена на $\widehat {S}$ аналитическими полиномами.

Поликруговая алгебра обладает этим свойством.
 Следующая теорема показывает, в частности,
что и  в случае, когда
полугруппа $S$  архимедово и линейно упорядочивает группу $G$  (что
равносильно тому, что $G$ есть подгруппа аддитивной группы $\Bbb{R}$ и
$S=G\cap\Bbb{R}_+$; см. \cite{Gam}),  алгебра $A(\widehat {S})$ тоже
обладает  свойством полиномиальной аппроксимации. В соответствии с
\cite{Mir} полугруппу $S$ будем называть конусом в $G$, если для любого
$x\in G$ найдется такой $\rho\in\widehat {S}_+$, что  $\tilde{\rho}(x)>1.$

{\bf Теорема 3.} {\it  Предположим, что $S$ есть конус в $G$, что
  $S^{-1}\cap S=\{e\}$, и что  полугруппа $S$
не содержит нетривиальных  простых идеалов. Тогда алгебра $A(\widehat {S})$ обладает
свойством полиномиальной аппроксимации}.

{\bf Доказательство}. Заметим сначала, что множество  $A(\widehat {S})$
содержится в пространстве $H^2(\widehat {S}\setminus X)$, определенном в
\cite{Mir}, т. е. что каждая функция  $F\in A(\widehat {S})$  обладает
следующими свойствами:

1) для каждого $\rho\in \widehat {S}_+\setminus X$ функция $F_{\rho}:\chi\mapsto
F(\rho\chi)$  принадлежит $L^2(X)$ и  нормы всех таких функций в $L^2(X)$
ограничены в совокупности;

2) для любых   $\rho, \rho_1\in \widehat {S}_+\setminus X$
преобразования Фурье ${\cal F}F_{\rho_1}$ и  ${\cal F}F_{\rho_1\rho^0}$
совпадают на группе частных полугруппы $S(\rho)$;

3) $F$ обобщенная аналитическая в смысле определения 1.

В самом деле, в доказательстве нуждается лишь свойство 2, которое
достаточно проверить для единственного полухарактера $\rho=\omega$,
принимающего нулевые значения. Но в этом случае оно сразу следует из
теоремы 2, примененной к функции   $\psi\mapsto  F(\rho_1\psi)$.

Если теперь мы положим $F^*=F|X$, то  следствие 5.2 из \cite{Mir}
показывает, что спектр (т. е. носитель преобразования Фурье) функции
$F^*$  содержится в $S$. Поэтому для любого $\varepsilon > 0$ найдется аналитический
полином $p$ такой, что
$$
|p(\chi)-F(\chi)|<\varepsilon \mbox{ для всех } \chi\in X.
$$
\noindent
Пусть $\Phi(\psi)=p(\psi)-F(\psi)\quad(\psi\in\widehat {S})$, и предположим, что
$\max_{\widehat {S}}|\Phi|=|\Phi(\psi_1)|$, где $ \psi_1=\rho_1\chi_1\in
\widehat {S}\setminus X$.  У нас $\rho_1 \ne 1$; кроме того, мы можем считать, что
$\rho_1 \ne \omega$, поскольку  в силу теоремы 2 $|\Phi(\omega)|\leq\max_X|\Phi|$.
 Обозначим через $\kappa$ отображение
множества $\Pi\cup\{0\}$ в  $\widehat {S}$, заданное формулой
$\kappa(z)=\rho_1^z\chi_1$. Тогда модуль аналитической в $\Pi$ функции
$\Phi\circ\kappa$ достигает своего максимума в точке $z=1$, а потому
$\Phi\circ\kappa=const$ в $\Pi$. С учетом непрерывности получаем
$\Phi\circ\kappa(0)=\Phi\circ\kappa(1)$, т. е.
$\Phi(\psi_1)=\Phi(\chi_1)$, так как у нас $\rho_1^0=1$. Таким образом,
$$
|p(\psi)-F(\psi)|<\varepsilon \mbox{ для всех } \psi\in\widehat {S}.\Box
$$

{\bf Теорема 4}. {\it Пусть алгебра $A(\widehat {S})$ обладает свойством
полиномиальной аппроксимации. Тогда  пространство максимальных идеалов
$M_{A(\widehat {S})}$  этой алгебры  гомеоморфно $\widehat
{S}$, а ее граница Шилова  $\partial_{A(\widehat {S})}$ гомеоморфна} $X$.

{\bf Доказательство}.  Для  $\varphi\in M_{A(\widehat {S})}$
положим $\psi(a):= \varphi(\widehat {a}), a\in S$. Тогда
$\psi \in \widehat {S}$, причем $\widehat{a}(\psi)= \varphi(\widehat {a})$.
 Если  $p=\sum_{i=1}^{n}c_i \widehat {a_{i}}$ ($c_{i}
\in \Bbb{C},a_{i}\in S$) --- аналитический полином, то
$$
\varphi(p)=\sum_{i=1}^{n}c_i \varphi(\widehat {a_{i}})=
\sum_{i=1}^{n}c_i \widehat {a_{i}}(\psi)=p(\psi).
$$
Для произвольной функции $F$ из $A(\widehat {S})$ имеем
$F=\lim_{n \to \infty} p_{n}$, где $ p_{n}$ --- аналитические полиномы.
Поэтому
$$
\varphi(F)=
\lim_{n \to \infty} \varphi(p_{n})=\lim_{n \to \infty}
p_{n}(\psi)=F(\psi).
$$

Рассмотрим отображение $\Lambda: \widehat {S} \longrightarrow
M_{A(\widehat {S})}$, которое каждому $\psi$ из $\widehat {S}$ сопоставляет
комплексный гомоморфизм  $\varphi_{\psi}$ из  $M_{A(\widehat {S})}$ по формуле
$\varphi_{\psi}(F)=F(\psi)$.
Тогда $\Lambda$ --- гомеоморфизм. Действительно, если
$\psi_{1}, \psi_{2}\in \widehat {S},\psi_{1} \not= \psi_{2}$,
то  $\varphi_{\psi_{1}}(\widehat {a}) \not=
\varphi_{\psi_{2}}(\widehat {a})$ для некоторого $a\in S$, и
отображение $\Lambda$  является сюрьективным по доказанному выше.
Таким образом, оно биективно. Для доказательства его непрерывности
рассмотрим последовательность $\psi_{n}\in \widehat {S}$, сходящуюся к
$\psi\in \widehat {S}$. Тогда
$\Lambda(\psi_{n})$ сходится к $\Lambda(\psi)$, так  как   при
 $ F\in A(\widehat {S})$ имеем
$\varphi_{\psi_{n}}(F) \to \varphi_{\psi}(F)$   в силу непрерывности $F$.
 Поскольку $\widehat {S}$ --- компакт, то $\Lambda$ ---
 гомеоморфизм.

  Докажем второе утверждение теоремы. Покажем сначала, что $\partial_{A(\widehat {S})}
\subset X$. Для этого достаточно доказать, что $X$ является границей. Если это не так, то

$$
M:=\max\limits_{\widehat {S}}|F|>m:= \max\limits_X |F|
$$
\noindent
для некоторой функции  $ F\in A(\widehat {S})$.  Для $\varepsilon<(M-m)/2$ подберем
аналитический полином $p$ таким образом, чтобы $\max_{\widehat
{S}}|F-p|<\varepsilon$.  Тогда $|p(\chi)|<|F(\chi)|+\varepsilon\leq m+\varepsilon$
при всех $\chi\in X$.
Так как любой комплексный гомоморфизм из границы Шилова алгебры   $l_1(S)$
продолжается до комплексного гомоморфизма содержащей ее алгебры $l_1(G)$,
то
$X$  содержит  границу Шилова алгебры $l_1(S)$.  А так как для любого
$\zeta\in \widehat {S}$ отображение $f \mapsto \widehat {f}(\zeta)$
есть комплексный гомоморфизм алгебры  $l_1(S)$, то
 $M_{l_1(S)}\supseteq \widehat {S}$ (на самом деле, имеет место равенство).
Поэтому    $\max_{\widehat {S}}|p|=
\max_X |p|$ (аналитические полиномы являются преобразованиями  Гельфанда функций из
$l_1(S)$ с конечными носителями). Таким образом,  $\max_{\widehat S}|p|\leq m+\varepsilon$.

С другой стороны, $|p(\psi)|>|F(\psi)|-\varepsilon$ при всех  $\psi\in
\widehat {S}$, а потому  $\max_{\widehat S}|p|\geq M-\varepsilon$ , что
противоречит выбору $\varepsilon$. Это доказывает требуемое включение.

Наконец, так как  алгебра $A(\widehat {S})$ инвариантна относительно
естественного действия группы $X$ (умножения на характеры являются автоморфизмами
полугруппы   $\widehat {S}$), то такова и ее граница Шилова, а потому
эта граница совпадает с $X$. $\Box$

{\bf Замечание.} Утверждение теоремы может быть выведено также из результатов таботы
\cite{AS}.

Ниже  $A(X)$ будет обозначать равномерную алгебру, состоящую из сужений на $X$ функций из
$A(\widehat {S})$.

{\bf Следствие}.  Пусть алгебра $A(\widehat {S})$ обладает свойством
полиномиальной аппроксимации. Тогда она изометрически изоморфна алгебре $A(X)$.

Действительно, требуемым изоморфизмом будет отображение сужения на $X.\Box$

\vspace{5mm}
\centerline{\bf 3 Интерполяционные множества для алгебры $A(\widehat {S})$}
\vspace{5mm}

Далее через ${\cal M}(X)$ будем обозначать алгебру всех регулярных комплексных борелевских
мер на компакте $X$. Для меры $\nu\in {\cal M}(X)$ запись $\nu\bot A(\widehat
{S})$ будет означать, что $\int Fd\nu=0$ для всех $F\in A(\widehat {S})$.
Через ${\cal M}_\omega$ обозначим множество всех вероятностных мер из ${\cal M}(X)$,
представляющих функционал $F\mapsto F(\omega), \quad F\in A(\widehat
{S})$. Теорема 2 в точности означает, что мера Хаара $\sigma$ группы $X$
принадлежит ${\cal M}_\omega$.

По аналогии со случаем шара (\cite{Rud84}, глава 10) введем шесть символов $(I), (PI), (Z),
(P), (N), (TN)$ для обозначения свойств,
которыми может обладать компакт $K \subset X$ по отношению к алгебре
$A(\widehat {S})$.

 $K$ будем называть $(I)$-{\it множеством
(интерполяционным множеством)}, если  любая комплекснозначная непрерывная
функция
на $K$ допускает продолжение до функции из $A(\widehat {S})$.

 $K$ будем называть $(PI)$-{\it множеством
(пик-интерполяционным множеством)}, если для любого $g \in C(K), g\not\equiv 0$,
 существует функция $ h \in A(\widehat {S})$  такая, что $g(\chi)~=~h(\chi)$
для всякого    $\chi \in K$,
  а для любого $\psi \in \widehat {S} \setminus K$ выполняется неравенство $|h(\psi)|<
\|g\|_{K}$, где $\|g\|_{K}=\max_K|g|$.

$K$ будем называть $(Z)$-{\it множеством
(нулевым множеством)}, если существует $f \in A(\widehat {S})$, такая, что
$f(\chi)=0$ для любого $\chi \in K$, и $f(\psi)\not= 0$ для любого $\psi \in \widehat {S}
\setminus K$.

$K$ будем называть $(P)$-{\it множеством
(множеством пика)}, если найдется функция $h \in A(\widehat {S})$   такая, что
$h(\chi)=1$
для любого $ \chi \in K$, и $|h(\psi)|~<~1$ для любого $ \psi \in \widehat {S} \setminus K$
($h$ --- пик-функция для $K$).

$K$ будем называть $(N)$-{\it множеством
(нуль-множеством для любой меры $\nu$ на $X$, аннулирующей $A(\widehat {S})$)},
 если равенство $|\nu|(K)=0$ выполняется для любой
меры $\nu\in {\cal M}(X)$ такой, что $\nu \bot A(\widehat {S})$.

 $K$ будем называть $(TN)$-{\it множеством (вполне нулевым множеством),}
 если равенство $\mu(K)=0$ выполняется для любой  представляющей
меры $\mu\in {\cal M}_\omega$.

Аналогичные понятия (кроме свойства $(TN)$) можно ввести и для алгебры $A(X)$.

 {\bf Теорема 5}. {\it Пусть алгебра $A(\widehat {S})$ обладает свойством
полиномиальной аппроксимации. Для любого компакта $K \subset X$ типа
$G_\delta$ такого, что дополнение $\widehat {S}\setminus K$ односвязно, свойства $(I),
(PI), (Z), (P), (N)$ равносильны и влекут свойство} $(TN)$.

{\bf Доказательство}.  Равносильность будем доказывать по схеме
$$
(I)\Rightarrow(N)\Rightarrow(PI)\Rightarrow(Z)\Rightarrow(P)\Rightarrow(I).
$$

$(I)\Rightarrow(N)$. Если $S$ счетно, то любая точка $\chi\in X$ будет точкой пика для
$A(\widehat {S})$ (а потому и для $A(X)$), что показывает пик-функция
$$
\psi\mapsto 1/2(1+\langle\psi,\chi\rangle),
$$
\noindent
построенная для веса $w>0$.  Так как $(I)$-множество $K$ алгебры
$A(\widehat {S})$  будет $(I)$-множеством и для алгебры $A(X)$,то требуемая импликация
сразу следует из теоремы Варопулоса (\cite{Rud84}, теорема 10.2.2), примененной к
алгебре $A(X)$. В случае произвольного $S$ можно воспользоваться методом
из \cite{Rud74}, с. 114-115.

$(N)\Rightarrow(PI)$.  В силу следствия теоремы 4 отображение сужения на $X$
есть изометрический изоморфизм алгебр $A(\widehat {S})$  и $A(X)$, а
потому $A(X)$ замкнута в $C(X)$. Теорема Бишопа (\cite{Rud84}, теорема
10.3.1), примененная к алгебре $A(X)$, показывает теперь, что каждое
$(N)$-множество $K$ будет $(PI)$-множеством для алгебры $A(X)$, а следовательно и
для  $A(\widehat {S})$.

$(PI)\Rightarrow(Z)$. Это очевидно.

$(Z)\Rightarrow(P)$. Пусть $F$ --- функция из  $A(\widehat {S})$ с
множеством нулей $K$, причем $\max_{\widehat {S}} |F|<1$. Так как
$\widehat {S}\setminus K$ односвязно, то $F$ обладает непрерывным
логарифмом $g$ на этом множестве, который, очевидно, аналитичен на
$\widehat {S}\setminus X$.  Кроме того, ${\rm Re} g<0$ и ${\rm Re}
g(\psi)\to -\infty$,   когда $\psi$ стремится к некоторой точке из $K$.
Следовательно, функция $h=g/(g-1)$ на   $\widehat {S}\setminus K$  и $h=1$
  на $K$  принадлежит $A(\widehat {S})$ и имеет пик на $K$.

$(P)\Rightarrow(I).$  Пусть  $R:A(\widehat {S})\to C(K)$ --- отображение
сужения на $K, A_K=ImR$. Нужно доказать, что $A_K=C(K)$, если  $K$  --- $(P)$-множество.

Если $J$ --- идеал функций из $A(\widehat {S})$, равных нулю на  $K$, то
фактор алгебра $A(\widehat {S})/J$ изометрически изоморфна $A_K$
(требуемым изоморфизмом служит $F+J\mapsto F|K$; см., напр.,\cite{Burb},
c. 120), а потому $A_K$ замкнуто в $C(K)$.

Покажем теперь, что $M_{A_K}$ можно отождествить с $K$.  Если $\varphi\in
M_{A_K}$, то $\varphi R\in M_{A(\widehat {S})}$, и, как показано в
доказательстве теоремы 4, существует такой $\psi\in \widehat {S}$, что
$(\varphi R)(F)=F(\psi)$ для любого $F\in A(\widehat {S})$. Если $h$ --- пик-функция для
$K$, то $Rh$ есть единица алгебры $A_K$, а потому  $h(\psi)=(\varphi
R)(h)=1$. Значит, $\psi\in K$. Если $g\in A_K$,то  $g=RF$ для некоторого
  $F\in A(\widehat {S})$, и, следовательно,
$$
\varphi(g)=(\varphi R)(F)=F(\psi)=g(\psi),
$$
\noindent
причем $\psi$ с таким свойством единственно, так как $A_K$ разделяет точки
 $K$. Отождествляя $\varphi$ с $\psi$, получаем равенство $M_{A_K}=K$.

Для $a\in S$ положим теперь  $\widehat{a}_K:= \widehat{a}|K$.  Поскольку элементы
  $\widehat{a}_K$ не обращаются в нуль на $M_{A_K}=K$
($|\widehat{a}_K|=1$),  то они обратимы в алгебре  $A_K$, причем
$(\widehat{a}_K)^{-1}=\overline{\widehat{a}_K}$. Поэтому подалгебра в
$A_K$, порожденная элементами
$\widehat{a}_K,\overline{\widehat{a}_K},\quad a\in S$, симметрична. Так как
она, очевидно, разделяет точки $K$, то $A_K=C(K)$ по теореме
Вейерштрасса-Стоуна.

Покажем, наконец, что  $(P)\Rightarrow(TN).$  Действительно, если $h$ --- пик-функция
для $K$, то   для любого натурального $n$ и для любой меры $\mu\in {\cal
M}_\omega$  имеем
$$
h^n(\omega)=\int_X h^n(\chi)d\mu(\chi).
$$
\noindent
Полагая здесь $n\to\infty$, получаем $0=\mu(K).\Box$

Теорема 5 обобщает известные результаты для поликруга (\cite{Rud74},
теорема 6.1.2).
Как и в \cite{Rud74}, получаем также

{\bf Следствие 1.}   Любое компактное $G_\delta$-подмножество
$(Z)$-множества в $X$ будет  $(Z)$-множеством в $X$.  Любое компактное
$G_\delta$-подмножество
множества  $X$, являющееся счетным объединением $(PI)$-множеств,  будет  $(PI)$-
множеством в $X$.

Действительно, оба свойства очевидны для $(N)$-множеств.$\Box$

{\bf Следствие 2.} $\sigma(K)=0$  для любого $(I)$-множества  $K$ ($\sigma$ ---
мера Хаара группы $X$).

Это следует из включения $\sigma\in {\cal M}_\omega$ (теорема 2). $\Box$

В \cite{AS} доказано, что для любого $\zeta\in\widehat {S}$  существует такая регулярная
борелевская мера $m_{\zeta}$ на $X$,  что
справедливо следующее обобщение формулы Пуассона
$$
\widehat {f}(\zeta)=\int\limits_X\widehat {f}dm_{\zeta}\quad (f\in l_1(S)).
$$
\noindent
При этом имеет место

{\bf Следствие 3.} $m_{\zeta}(K)=0$  для любого $(I)$-множества  $K$, если
$\zeta\in \widehat {S}\setminus K$.

{\bf Доказательство}. Пусть $K$ есть множество нулей функции $F\in A(\widehat
{S})$. Если $m_{\zeta}(K)>0$, то
$$
\int_{X} \log|F|dm_{\zeta}=-\infty,
$$
\noindent
и в силу теоремы 3.5 из \cite{Ar} $F(\zeta)=0$, что противоречит выбору $\zeta.\Box$

С помощью теоремы 5 можно описать $(I)$-множества в $X$  в случае, когда
$S$ линейно и архимедово упорядочивает группу $G$. При этом не нарушая
общности мы будем считать, что $G$ есть подгруппа аддитивной группы
$\Bbb{R}, \quad S=G\cap \Bbb{R}_+$.
Тогда $\Bbb{R}$  изоморфна некоторой подгруппе группы $X$, которую мы тоже
обозначим $\Bbb{R}$ ( см. \cite{Gam}, глава 7, $\S 4$). В этой ситуации групповую
операцию в $X$ принято  записывать аддитивно.

{\bf Следствие 4.} (Ср. \cite{Gam}, с. 254-255).  Пусть $S$ линейно и архимедово
упорядочивает группу $G$.    Компактное $G_\delta$-подмножество
$K$ множества  $X$ будет  $(I)$-множеством  тогда и только тогда, когда
для любого $\chi\in X$  пересечение $K\cap(\chi+\Bbb{R})$ имеет нулевую линейную меру.

{\bf Доказательство}. Алгебра $A(\widehat {S})$ обладает свойством
полиномиальной аппроксимации, например,  по теореме 3 (условия этой теоремы
выполняются в силу свойств полухарактеров полугруппы $S$, установленных
в \cite{Hof}, c. 449).Тогда $M_{A(\widehat {S})}=\widehat{S}$.
Покажем, что для любого подмножества $K\subset X$
множество  $\widehat {S}\setminus K$ односвязно. Известно  (\cite{Gam}, теорема 7.4.1), что
$\widehat {S}$ можно отождествить с пространством
$$
conX:=X\times[0,1]/X\times\{0\},
$$
\noindent
и при этом группа $X$ отождествляется с $X\times\{1\}$. Точку из $conX$, отвечающую
$X\times\{0\}$, обозначим $\omega$. Заметим, что тождественное отображение
множества $conX\setminus K$ гомотопно постоянному отображению этого множества в $\omega$.
Требуемую гомотопию дает отображение $f_t:conX\to conX$, определяемое
формулой $f_t(\chi,r)=(\chi,tr)$ при $t,r\in [0,1], \chi\in X$. Таким
образом, множество $conX\setminus K$ стягиваемо, а потому односвязно. В
силу теоремы 5 свойства $(I)$ и $(N)$ равносильны, и осталось
воспользоваться следствием 7.5.2 из \cite{Gam}. $\Box$

{\bf Следствие 5}. Пусть замкнутая подгруппа $H\subset X$ такова, что
группа $X/H$ метризуема. Тогда $H$ будет $(I)$-множеством если и только
если множество $H^{\bot}\cap S$ (рассматриваемое как множество функций на
$X/H$) разделяет точки $X/H$.

Это  следует из теоремы 5 и теоремы 6 в \cite{GT2}.

Теорема 5  позволяет также установить аналог теоремы Деви-Эксендала (см.
\cite{Rud84}, теорема 10.4.3), дающий достаточный признак
интерполяционности множества в терминах его покрытий. Для $\zeta\in X,
\delta>0$ и веса $w>0$ положим
$$
V(\zeta,\delta)=\{\chi\in X: |1-\langle\chi,\zeta\rangle|<\delta\}
$$
\noindent
(угловые скобки обозначают, как и раньше, скалярное произведение в
$l_2(S,w)$).

{\bf Теорема 6}. {\it В условиях теоремы 5 предположим, что $S$ счетно и что
 компакт $K$  обладает также следующим свойством:
\par
для любого $\varepsilon >0$ существует конечное число множеств
$V(\zeta_i,\delta_i),\quad i=1,\ldots, m$ с $\zeta_i\in K$ и
$\sum\delta_i<\varepsilon$, таких, что
$$
K\subset V(\zeta_1,\delta_1)\cup\ldots\cup V(\zeta_m,\delta_m).
$$
Тогда $K$ есть  $(PI)$-множество}.

Действительно, по-существу повторяя  доказательство теоремы Деви-Эксендала, получаем, что
$K$ является  $(N)$-множеством. $\Box$

При $S=\Bbb{Z}_+^n$ теорема 6 есть ослабленный вариант теоремы Форелли
(\cite{Rud74}, c. 126, следствие).

Укажем также следующее применение теоремы 6. Ниже $I=[0,1],\quad w$ ---
фиксированный положительный вес на счетной полугруппе $S$.

{\bf Определение 3.}  Отображение $\gamma:I\to l_2(S,w)$ класса $C^1$
будем называть {\it комплексно-касательной кривой на } $X$, если
$\gamma(I)\subset X$, и при всех  $t\in I$

$$
\langle\gamma(t),\gamma ^{'}(t)\rangle=0.
$$

{\bf Теорема 7}. {\it Пусть $S$ счетно. В условиях теоремы 5 образ $K=\gamma(I)$
комплексно-касательной  кривой $\gamma$ на $X$  будет  $(PI)$-множеством}.

{\bf Доказательство}. Достаточно убедиться, что $\gamma(I)$ удовлетворяет
условию теоремы 6. Пусть $M>\max_{s,t\in I}|\langle\gamma^{'}(s),\gamma ^{'}(t)\rangle|$.
Для любого $\varepsilon >0$ рассмотрим разбиение $\{a_i\}$ отрезка $I$ на
$m>M/\varepsilon$ равных  частей, и для любого $t\in I$ выберем $a_i\leq t$ так,
что $t-a_i\leq 1/m$. Теорема будет доказана, если мы установим, что
$\gamma(t)\in V(\gamma(a_i),\delta)$, где $\delta=M/m^2$. Но это вытекает
из следующей оценки
$$
|1-\langle\gamma(t),\gamma(a_i)\rangle|=|\langle\gamma(t),\gamma(t)-\gamma(a_i)\rangle|=
|\langle\gamma(t),\int\limits_{a_i} ^t\gamma ^{'}(s)ds\rangle| =
$$
$$
|\int\limits_{a_i} ^t\langle\gamma(t),\gamma ^{'}(s)\rangle ds| =
|\int\limits_{a_i} ^t\langle\gamma(t)-\gamma(s),\gamma ^{'}(s)\rangle ds|=
$$
$$
|\int\limits_{a_i} ^t\langle\int\limits_s^t\gamma^{'}(r)dr,\gamma ^{'}(s)\rangle ds| =
|\int\limits_{a_i} ^tds\int\limits_s^t\langle\gamma^{'}(r),\gamma
^{'}(s)\rangle dr| \leq
$$
$$
\int\limits_{a_i} ^tds\int\limits_s^t|\langle\gamma^{'}(r),\gamma
^{'}(s)\rangle| dr <  M(t-a_i)^2\leq\delta.\Box
$$

Отметим, что  не всякая кривая на $X$ класса $C^1$ является
$(I)$-множеством, в чем нас убеждает следующий

{\bf Пример 2}. Пусть существует такой полухарактер $\rho\in \widehat
{S}_+$, что  $\rho(s)>0$ при всех $s\in S$. Тогда отображение $\gamma:I\to
l_2(S,w)$, заданное формулой $\gamma(t)=\rho^{it}$, есть кривая класса
$C^1$ на $X$, но $\gamma(I)$ не является $(Z)$-множеством. Действительно, пусть
$F\in A(\widehat {S}),\quad F|\gamma(I)=0$. Тогда функция $f(z):=F(\rho^z)$ будет
  аналитической в открытой правой полуплоскости,  непрерывной
в ее  замыкании и равной нулю на отрезке $iI$. Следовательно, $f=0$ по
граничной теореме единственности.

\vspace{35mm}

\label{liter}\def\refname{}

\centerline{{\bf Литература}}

{\it Гомельский государственный университет\\
им. Ф. Скорины \ (Беларусь)}

\end{document}